\documentclass[12pt]{article}
\usepackage{color}
\usepackage{amssymb,latexsym}
\usepackage{pdfsync}
\usepackage[colorlinks,backref=page
]{hyperref}
\usepackage{amsmath,amsthm}
\usepackage{indentfirst}
\theoremstyle{plain}
\usepackage{amsmath,amsthm}
\usepackage[colorlinks, backref=page]{hyperref}
\usepackage{cleveref}
\theoremstyle{plain}

\newtheorem{thm}{Theorem}[section]

\newtheorem{lem}{Lemma}

\theoremstyle{definition}
\newtheorem{defi}{Definition}

\newtheorem{rem}{Remark}
\numberwithin{equation}{section}
\allowdisplaybreaks

\def \d {\mathrm{d}}
\def \Ric {\mathrm{Ric}}
\def \vol{\mathrm{Vol}}

\title{\Large \bf \boldmath\ \\ Liouville Type Theorem about $p$-harmonic 1 form,  $p$-harmonic map and  harmonic $ q $ form } 

\author{\large  Xiangzhi Cao \thanks{School of information engineering, Nanjing Xiaozhuang University, Nanjing 211171, China}\thanks{aaa7756kijlp@163.com}} 

\date{}

\begin{document}
	
	\maketitle
	\tableofcontents
	
	\begin{abstract}
		In this paper, we will use the 
			normalized intetral Ricci curvature to investigate Liouville type property of $ p $ harmonic function on Riemannian manifold. secondly, we will use the  BiRic curvature to obtian Liuville theorem for $ p $ harmonic function or $ p $ harmonic 1 form.  Lastly, we  we will use the  BiRic curvature to obtian Liouville theorem forharmonic $ q(q\geq 2)$ form.
	\end{abstract}
	
	\section{Introduction}
	Let $ (M,g) $ and $ (N,h) $ be two Riemannian manifold, $ u: (M,g) \to (N,h) $
	The energy functional of $p$ harmonic map is  defined by
	\begin{equation*}
		E_{p}(u)=\int_{M}\frac{| \nabla u|^{p}}{2}\mathrm{d} v_g,
	\end{equation*}
	whose Euler-lagrange equation is as follows:
	\begin{equation}\label{bitensionfild}
		\tau_{p}(u)=\mathrm{div}(|du|^{p-2}du),
	\end{equation}
	In particular, when $ N=\mathbb{R}, $  $ p $ harmonic map reduces to $ p $ harmonic function.
	It is well known that the regularity of  $ p $ harmonic function is not better than $  C^{1,\alpha} $ (c.f.\cite{wang2011local}\cite{Tolksdorf1984}\cite{lindqvist2017notes} ).
	
	There are a lot of works on Liouville theorem for $ p $ harmonic map and p harmonic function. Generally speaking, the common strategy is to use Bochner formula and Kato inequality  by imposing the curvatre of the domain manifold and target manifold. 	In \cite{wangql}, Wang  studied p harmonic map and p harmonic funtion from submanifold in partially negative manifold to Cartan Hardmard manifold. In \cite{MR3692378}, we extended the results of  \cite{wangql} to $ p $ harmionic function and $ p $ harmonic map.   One can   also refer to \cite{p-harmonicmap,Kawai1999,prs,MR1145657,Zhang2016A,Baird1992,Li2001,dung2017p,pigola2008constancy} and reference therein for the works on Liouville theorem of   $ p  $ harmonic maps.
	
	In previous literatures,  Liouville type theorem for $ p $-harmonic map was obtained under certain conditions, for example, lower bound on Ricci curvature of $ M $ ,   the Sobolev inequality or weighted Poincare inquality holds on $ M $. If $ M $ is  a submanfold in Riemannian  manifold, then  the second fundermentan form shall satisfy some conditons, meanwhile $ M $ is stable submanfold. However, little is known for Liouville  type theorem of $ p $ harmonic map if the domain manifold has small normalized integral Ricci curvature(c.f. \eqref{norm}).

	Integral Ricci curvature is an important notion which can be used to study comparison theom(\cite{zbMATH07024088}\cite{zbMATH01123720}) and gradient estimates. In \cite{zbMATH02114455}, Dai used the integral Ricci curvature to study the heat kernel bound.
	In \cite{zbMATH07178581}, Wang and Wei used the notion of integral Baky-\'Emery Ricci curvature to obtain the local sobolev inequality.

	In this paper, we will prove Liouville type theorem for $  p $ harmonic function on manifold M whose normalized integral Ricci curvature is sufficiently small.
		
	Let $ (M,g, e^{-f}\mathrm{d}v_g) $ be metric measure space, let $ \pi: E\to M $ be a vector bundle of rank m over$  M $ , we call $ \omega $ is $ E $-valued weighted $ p $ harmonic form if it satisfies 
	
	$$ d \omega=0, \delta_f(|\omega|^{p-2}\omega)=0.$$
	where $ \delta_{f}=\delta +i_{\nabla f}. $ One can refer to \cite[section 3]{zbMATH07180899} for the definition of the two operators $ d $ and $ \delta $. In fact, they are just a slightly generaliztions of exterior differential operator and  codifferential operator on Riemannian manifold. 
	
	When $ f=0 $,		we call $ \omega $  $ E $ valued $ p $ harmonic form if it satisfies	
	$$ d \omega=0, \delta(|\omega|^{p-2}\omega)=0.$$
	Moreover, let $ (M,g) \to (N,h) $, take $ E= u^{-1}TN ,$	then it is easy to see that $ du $ is $ u^{-1}TN $ valued $ p $ harmonic form on $ M $.
	
	It is interesting to obtain Liouville theorem for p harmonic form  on Riemannian manifold.	The research of $ p $ harmonic $ l $ form mainly focus on vanishing property and finiteness of the vector space  $ p $ harmonic $ l $ form. On the one hand, there are a lot researches  for vanishing theorem of $ p $ harmonic $ l $ form on submanifold in Riemannian  manifold with special curvature property, such as \cite{2}\cite{zbMATH06859577}  \cite{zbMATH07040685}. On the other hand, the present research mainly focus on the vanishing theorem of $ p $ harmonic $ l $ form on manifold with special structure such as conformally flat or special curvature property, such as \cite{zbMATH06827101}\cite{9}\cite{zbMATH06668693}\cite{zbMATH07361892}.
	One can refer to \cite{2} for the studies on finiteness property of vector space p harmonic  l form on Riemannian manifold.

	Seo and Yun  \cite{zbMATH07180899} proved that  if  manifold has nongegatrive Bakry-\'Emery Ricci curvature, then weighted $ p $ harmonic $ 1 $ form must vanish. As an corollary, they got the Liouville theorem of weighted p harmonic map on manifold with nonnegative Bakry-\'Emergy Ricci curvature .

	In \cite{shen1996stable} \cite{shen1997geometry}, Shen and Ye introduced the concept of BiRic  curvature. Bi-Ricci curvature is an important  tool to obtain vanishing theorem for diferential forms and to estimate Ricci curvature of submanifold when the Bi-Ricci curvature is bounded below by a suitable constant. Bi-Ricci curvature can also be used to study stable  submanifolds immersed in a Riemannian manifold. Moreover, nonnegative BiRic curvature is weaker than  nonnegative sectional  curvature.

	Tanano \cite{tanno1996l2} obtained Liouville theorem of $ L^2 $ harmonic 1-forms on $ M $ which complete noncompact orientable stable minimal hypersurface in a Riemannian manifold with non-negative BiRic curvature.	Wang \cite{wang2006harmonic}  proved the Liouville theroem for harmonic map with finite energy in terms of BiRic curvature.	Cheng \cite{cheng20002} used the condtion that  $ BiRic > \frac{n-5}{4}H^2 $ to improve the results in \cite{zbMATH01208191} , and obtained the vanishing theroem .
	Li and Zou \cite{zoulibiricci} further gerneralized the results in \cite{cheng20002} by relaxing the lower bound of BiRic. 
	Cheng \cite{cheng2006constant} studied the topologic property of   complete finite index hypersurface immersed in a manifold  $ N $ with constant mean curvature, if $ BiRic > \frac{n-5}{4}H^2 $.	We refer the readers to	 the just mentioned references and references therein for further discussions.
	
	One  can refer to \cite{zbMATH07153791} for the  tensor.
	\begin{equation*}
		\begin{split}
			BiRic_f^{\delta}(X,n):=	Ric_f(X,X)+\delta Ric_f(n)-K(X,n).
		\end{split}
	\end{equation*}
	where $ X,n $ are two unit vector field. When $ f=0,\delta=1, $ it is just the BiRicci tensor.

	Let $M$ be an $n(\geq 3)$-dimensional complete $ f $-stable  noncompact $ f $-minimal submanifold isometrically immersed in an $(n + k)$-dimensional metric measure space $(\bar{M},g, e^{-f}\mathrm{d}v_g)$. In \cite{zbMATH07153791}, they proved that  if   $\bar{M}$    has $\overline{BiRic_f}\geq 0$  ,  any  wighted harmonic $1$-form with finite $L^{2}$-energy on $M$  vanishes.
	
	However, up to now, there is rare study about $ p $  harmonic map ($ p\neq 2 $) in terms of BiRic curvature or  weighted BiRic curvature.  In this paper, we will prove that  if  $\bar{M}$    has $\overline{BiRic}^{\delta}\geq 0$,  any  $p$-harmonic $1$-form with finite $L^{2p-2}$-energy on $M$  vanishes.

Vanishing property of harmonic  form is   important topic in differential geometry,  one can refer to these works (\cite{2011On,7,3,14}) and reference therein.  However, the condtions which were used to obtain vanishing theorem of harmonic form on submanifold or Riemmanian manifold   are, for example, requiring the lower bound of Weitzenböck curvature operator, or  requiring  Sobolev inequality holds on $ M $ or wighted poincar\'e inequality holds on $ M $, or imposing some  lower bound on the first eigenvalue of the operator $ \Delta_g $. If harmonic forms on  Riemqannian submanifold  are considered,  one may require the submanifold is stale or impose some conditions on the second fundermental form or mean curvature of submanifold. Of course, there are many extra conditions to obtian vanishing propertgy of harmonic form on manifold. Generally speaking, the studies on harmonic $ q(q\geq 2) $ form are more difficulty than that of $ L^2 $ harmonic $ 1 $ form. For harmonic 1 form in terms of BiRic curvature, on can refer to  \cite{cheng20002}   \cite{zbMATH01208191} and \cite{zbMATH07153791}. However, at present, there is rare study about  harmonic $ q(q\geq 2) $ form in terms of BiRic curvature or  weighted BiRic curvature.In this paper, we will prove that  if  $\bar{M}$    has $\overline{BiRic}^{\delta}\geq 0$,  any $ L^2 $ harmonic $ q $-form   on $M$  vanishes.

	The orgalization of this paper is as follows: in section 2, we give some Lemmas, formula and definitions in the proof of our Theorems. in section 3, we will derive Liouville theorem of p harmonic function on Riemannian manifold interms of integral curvature; in section 4, we will obtain Liouville theorem of p harmonic function or p harmonic 1 form in terms of BiRic curvature; in section 5 , we will study vanishing property for  harmonic $ q $ form in terms of BiRic curvature.
	
	\section{Preliminary}
	In this section, we give some lemmas which will be  used in this paper.
	We firstly recall  Kato inequality for p-harmonic function (cf. lemma 2.4 in \cite{zbMATH06451355}):
	\begin{equation}\label{2.1.2}
		|\nabla(|du|^{p-2}du) |^2\geq \frac{n}{n-1}|\nabla|du|^{p-1} |^2.
	\end{equation}
	However, for $ p $ harmonic map, we know that  (\cite{zbMATH04182289})	
	\begin{equation*}
		\begin{split}
			|\nabla(|du|^{p-2}du) |^2\geq |\nabla|du|^{p-1} |^2.
		\end{split}
	\end{equation*}

	A map $u$ is called  $L^{q}$-finite energy if $\int_{M}|\nabla u|^{q}\mathrm{d}v_g<\infty$.
	
	In the literature (c.f. \cite{zbMATH07024088}\cite{zbMATH07178581}) , we can find  the quantity
	\begin{equation*}
		\begin{split}
			Ric_{f-}^{H}=\max \{0,(n-1) H-\rho_f(x)\},
		\end{split}
	\end{equation*} where $\rho_f(x)$ is the smllest eigenvalue of $ Ric_f $ tensor of $M$.  
	\begin{rem}
		$\mathrm{Ric}_{f-}^{H}=0$ if and only if  $Ric_{f} \geq(n-1) H$.
	\end{rem}
	\begin{defi}[c.f. \cite{zbMATH07024088}]
		\[
		\|\phi\|_{p}(r):=\sup _{x \in M}\left(\int_{B(x, r)}|\phi|^{p} \cdot \mathcal{A}_{f} e^{-a t} \mathrm{d} t \mathrm{d} \theta_{n-1}\right)^{\frac{1}{p}},
		\]
		where $\partial_{r} f \geq-a$ for some constant $a \geq 0$, along a minimal geodesic segment from $x \in M$. Here $\mathcal{A}_{f}(t, \theta)$ is the volume element of weighted form $e^{-f} \mathrm{d} v_{g}=$ $\mathcal{A}_{f}(t, \theta) \mathrm{d} t \wedge \mathrm{d} \theta_{n-1}$ in polar coordinate, and $\mathrm{d} \theta_{n-1}$ is the volume element on unit sphere $S^{n-1}$. Sometimes it is convenient to work with the normalized curvature quantity(c.f. \cite{zbMATH07024088}\cite{zbMATH07178581})
	\begin{equation}\label{norm}
		\begin{split}
				\bar{k}(p,f, H, a, r):=\sup _{x \in M}\left(\frac{1}{V_{f}(x, r)} \cdot \int_{B(x, r)}\left(\operatorname{Ric}_{f}^{H}\right)^{p} \mathcal{A}_{f} e^{-a t} \mathrm{d} t \mathrm{d} \theta_{n-1}\right)^{\frac{1}{p}},
		\end{split}
	\end{equation}
		where $V_{f}(x, r):=\int_{B(x, r)} e^{-f} \mathrm{d} v$. The reader  should care about the difference of the same notation $ \bar{k}(p, H, a, r) $  in \cite{zbMATH07024088}\cite{zbMATH07178581}. Obviously, $\left\|\operatorname{Ric}_{f}^{H}\right\|_{p}(r)=0$ (or $\bar{k}(p, H, a, r)$ $=0$ ) iff $\operatorname{Ric}_{f} \geq(n-1) H$. When $f=0$ (and $a=0$ ), all above notations recover the usual integral curvature on manifolds.
	\end{defi}

	One can refer to these works \cite{zbMATH02114455} \cite{zbMATH01123720} for integral Ricci curvature.
	A $\operatorname{map} u$ is called  $ f $-weighted $L^{q}$-finite energy if $\int_{M}|\nabla u|^{q}e^{-f}\mathrm{d}v_g<\infty$.
		%

	\begin{lem}[c.f.  (2.5) in \cite{zbMATH06505276}]
		For anydifferential  $  q $ form, assume that  $ M $ is submanifold of $  N $, $ M $ has flat normal bundle and $ N $ has pure curvature tensor.
		\begin{equation*}
			\begin{split}
				& \sum\left( R_{i j}+f_{ij}\right) \omega^{ ii_{2} \cdots i_{q}} \omega^j_{ i_{2} \cdots i_{q}}
				-\frac{q-1}{2} \sum R_{k j i h} \omega_{i_{3} \cdots i_{q}}^{k{ }{j}} \omega^{i h{i_{3} \cdots i_{q}}}=F_{1}(\omega)+F_{2}(\omega),
			\end{split}
		\end{equation*}
		where
		\begin{equation}
			\begin{split}
				F_{1}(\omega)=& \sum_{k=1 i, i_{2}, \ldots, i_{q}}^{n} \bar{R}_{i k i k} \omega^{i i_{2} \ldots i_{q}} \omega_{i_{2} \ldots i_{q}}^{i} \\
				&-\frac{q-1}{2} \sum\left(\bar{R}_{i j i j} \omega^{i j i_{3} \ldots i_{q}} \omega_{i_{3} \ldots i_{q}}^{i j}+\bar{R}_{i j j i} \omega^{i j i_{3} \ldots i_{q}} \omega_{i_{3} \ldots i_{q}}^{j i}\right) \\
				=& \sum_{k=1}^{n} K_{i k} \omega^{i i_{2} \ldots i_{q}} \omega_{i i_{2} \ldots i_{q}}-(q-1) \sum K_{i j} \omega^{i j i_{3} \ldots i_{q}} \omega_{i j i_{3} \ldots i_{q}} \\
				=& \sum_{k=1}^{n} K_{i_{1} k} \omega^{i_{1} i_{2} \ldots i_{q}} \omega_{i_{1} i_{2} \ldots i_{q}} \\
				&-\sum_{n}^{n}\left(K_{i_{1} i_{2}}+K_{i_{1} i_{3}}+\cdots+K_{i_{1} i_{q}}\right) \omega^{i_{1} i_{2} i_{3} \ldots i_{q}} \omega_{i_{1} i_{2} i_{3} \ldots i_{q}} \\
				=& \sum_{h=q+1} K_{i_{1} i_{h}} \omega^{i_{1} i_{2} \ldots i_{q}} \omega_{i_{1} i_{2} \ldots i_{q}} \\
				=& \frac{1}{q} \sum_{t=1}^{q} \sum_{h=q+1}^{n} K_{i_{t} i_{h}} \omega^{i_{1} i_{2} \ldots i_{q}} \omega_{i_{1} i_{2} \ldots i_{q}} \\
				\geq & \frac{1}{q}\left(\inf _{i_{1}, \ldots, i_{n}} \sum_{t=1}^{q} \sum_{h=q+1}^{n} K_{i_{t} i_{h}}\right)|\omega|^{2} .	
			\end{split}
		\end{equation}
		and 
		\begin{equation*}
			\begin{split}
				F_2 \geq \frac{1}{2}(n^2|H |^2-C_{n,q}|A |^2)|\omega |^2,
			\end{split}
		\end{equation*}
		where $ C_{n,q}=\max\{q,n-q\} $
	\end{lem}
	
	\begin{lem}[c.f. Lemma 2.4 in \cite{zbMATH06859577}  or Lemma 2.5 in  \cite{MR3849353}]
		For any closed $q$-form $\omega$ and $\varphi \in \mathcal{C}^{\infty}(M)$, we have
		\[
		|d(\varphi \omega)|=|d \varphi \wedge \omega| \leq|d \varphi| \cdot|\omega| .
		\]
	\end{lem}	
	
	Let 
	\begin{equation*}
		\begin{split}
			A_{p,n,q}=\begin{cases}
				1+ \frac{1}{\max\{q,n-q\}}	& \text{if}\quad  p=2\\
				1+\frac{1}{(p-1)^2}\min\{1, \frac{(p-1)^2}{n-1}\}	&  \text{if} \quad  p>2 \quad  \text{and} \quad q=1\\
				1	& \text{if}\quad   p>2 \quad  \text{and} \quad 1< q\leq n-1
			\end{cases}
		\end{split}
	\end{equation*}
	
	\begin{lem}[Kato's inequality,cf. Lemma 2.2 in \cite{MR3849353} ]\label{kato}
		For $p \geq 2, q \geq 1$, let $\omega$ be an $p$-harmonic $q$-form on a complete Riemannian manifold $M^{n}$. The following inequality holds
		\[
		\left|\nabla\left(|\omega|^{p-2} \omega\right)\right|^{2} \geq\left.\left. A_{p, n, q}|\nabla| \omega\right|^{p-1}\right|^{2}.
		\]
		Moreover, when $p=2, q>1$ then the equality holds if and only if there exixts a 1-form $\alpha$ such that
		\[
		\nabla \omega=\alpha \otimes \omega-\frac{1}{\sqrt{q+1}} \theta_{1}(\alpha \wedge \omega)+\frac{1}{\sqrt{n+1-q}} \theta_{2}\left(i_{\alpha} \omega\right) .
		\]
		
	\end{lem}
	
	\section{ Liouville type theroem in terms of integral curvature condition}
	\begin{lem}[c.f. Propostion 1.3 in \cite{zbMATH07178581}]\label{269}
		Let $(M,g,e^{-f} \mathrm{d}v_g)$ be a complete smooth metric measure space. Assume that $\partial_{r} f \geq-a$ along all minimal geodesic segments for some constant $a \geq 0$. For $p>\frac{n}{2}$, there exists $\varepsilon=\varepsilon(n, p, a)>0$ such that if $r^2\bar{\kappa}(p, f, a, r) \leq \varepsilon$,
		Then the Sobolev inequality holds
		\[
		\int_{B_{x}(R)}|\nabla h| e^{-f} d \mathrm{vol} \geq 10^{-2 n} e^{-2 a} R^{-1}\left(\int_{B_{x}(R)} h^{\frac{n}{n-1}} e^{-f} d \mathrm{vol}\right)^{\frac{n-1}{n}},
		\]
		holds for all $h \in C_{0}^{\infty}\left(B_{x}(R)\right)$.
	\end{lem} 
	

	\begin{thm}
		Let $(M,g)$ be an $n(\geq 3)$-dimensional complete noncompact simply connected complete metric measure space. If there exists a postive constant $\epsilon$, such that  $ R^2\bar{k}(n, 0,H, 0, R)\leq \epsilon $ and    $\frac{1}{n-1}\geq (\frac{p-2}{p-1})^{2},p\geq 2$, $ 0<H $. We also assume that there is a point $q_1 \in M$ such that $Ric(X,X)|_{q_1}\neq 0 $ for all $ 0 \neq X \in T_{q_1}M . $  If one of the following conditions holds
		
		(1)$ H< \frac{1}{n-1} \bigg[\frac{1}{n-1}- (\frac{p-2}{p-1})^{2}\bigg] \frac{\epsilon_1\epsilon_2+1}{\epsilon_1^2} ,$ where  $\epsilon_1, \epsilon_2$ are  arbitrary small positve constants.
		
		(2) $ \vol(B(R)) \leq CR^{n-\epsilon} $, where $ C $ is a constant independent of $ R $, $\epsilon$ is arbitrary small positve constant.

		then  any $ p $-harmonic function with finite $L^{2 p-2}$-energy on $M$ is a constant.
		
	\end{thm}

	\begin{rem}
		$ R^2\bar{k}(n, 0,H, 0, R)$ is scaling invariant with regard to $ R $. $ R^2\bar{k}(n, 0,H, 0, R) \leq \epsilon$ doesn't imply that $ M $ is  compact. However, the condtion  $ \bar{k}(n, 0,H, 0, R) \leq \epsilon $ may imply $ M $ is compact.
	\end{rem}
	
	\begin{rem}
		The proof depends on the good bound in Kato's inequality, thus the theorem can not be applied to $ p $ harmonic map or $  p $ harmonic $ k $ form $ (k\geq 2) $. However, our methods can be applied to study $ p $ harmonic $ 1 $ form on manifold. The strategy is almost the same,  we omit it.
	\end{rem}
	\begin{proof}
		
		Here we modify the proof in \cite{MR3692378}. 
		First we recall  Bochner formula  for $p$-harmonic function(cf.lemma 2.1 in \cite{zbMATH06451355}):
		for any smooth function $u$, we have
		\begin{align} \label{2.1.1}
			&\frac{1}{2}\Delta|du|^{2p-2} =|\nabla|du|^{p-2}du|^{2}-\langle|du|^{2p-2}du,\Delta(|du|^{p-2}du)\rangle\nonumber\\
			&+|du|^{2p-4}\langle du(Ric^{M}(e_{k}),du(e_{k})\rangle ,
		\end{align}
		and the following  Kato inequality for $ p $-harmonic function (c.f. \cite[lemma 2.4]{zbMATH06451355} or \cite{MR3692378}):
		\begin{equation}\label{2.1.2}
			|\nabla u|^{2}\geq \frac{n}{n-1}|\nabla|u||^{2}.
		\end{equation}

		If we set $\lambda=|d u|^{p-1}$, we have(c.f. \cite{MR3692378})
		\begin{equation}\label{ev}
			\begin{split}
				\frac{1}{2} \Delta \lambda^{2}
				\geq& \frac{n}{n-1}|\nabla \lambda|^{2}-\left\langle|d u|^{2 p-2} d u, \delta d\left(|d u|^{p-2} d u\right)\right\rangle+|d u|^{2 p-4}\left\langle d u\left(\operatorname{Ric} c^{M}\left(e_{k}\right), d u\left(e_{k}\right)\right\rangle\right.\\
				\geq &\frac{n}{n-1}|\nabla \lambda|^{2}-\left\langle|d u|^{2 p-2} d u, \delta d\left(|d u|^{p-2} d u\right)\right\rangle+\left((n-1) H-\operatorname{Ric}_{-}^{H}\right) \lambda^{2}	
			\end{split}
		\end{equation}

		Now we choose a cut-off function on $M$  , which satisfies that (c.f. c.f. \cite{MR3692378} or \cite{wangql}):
		\begin{equation}\label{f9}
			\begin{cases}
				0\leq\phi(x)\leq1 ,& x\in M \\ \phi(x)=1, & x\in B_{R}(x_{0}) \\ \phi(x)=0 ,&x\in M-B_{3R}(x_{0})\\
				| \nabla \phi(x)|\leq \frac{C}{R}, & x \in M.
			\end{cases}
		\end{equation}
		Integrating on $M$ after multiplying both sides of $(2.3)$ by $\phi^{2}$, we get
		\begin{equation*}
			\begin{split}
				&\int_{M} \frac{1}{2} \phi^{2} \Delta \lambda^{2} \mathrm{d}v_g+\int_{M} R i c_{-}^{H} \lambda^{2} \phi^{2}\mathrm{d}v_g\\
				\geq &\frac{n}{n-1} \int_{M} \phi^{2}|\nabla \lambda|^{2} \mathrm{d}v_g-\int_{M} \phi^{2}\left\langle|d u|^{2 p-2} d u, \delta d\left(|d u|^{2 p-2} d u\right)\right\rangle \mathrm{d}v_g\\
				&+(n-1) H \int_{M} \lambda^{2} \phi^{2} \mathrm{d}v_g. 	
			\end{split}
		\end{equation*}

		It  is not hard to get that 
		\begin{equation*}
			\begin{split}
				&\frac{1}{Vol(B_{3R}(x))^{}}	\int_{M}  Ric _{-}^{H} \phi^{2} \lambda^{2} \mathrm{d}v_g \\
				&\leq  \frac{1}{Vol(B_{3R}(x))^{\frac{1}{n}}}\|  Ric_{-}^{H} \|_{L^{n}(B_{3R}(x))} \frac{1}{Vol(B_{3R}(x))^{\frac{n-1}{n}}}\| \phi^{2} \lambda^{2} \|_{L^{\frac{n}{n-1}}(B_{3R}(x))}\\
			\end{split}
		\end{equation*}	
		
		Thus, in Lemma \ref{269} , we take $ f=0 $ and $ a=0 $, the condition that  $ R^2\bar{k}(n, 0,H, 0, R)\leq \epsilon $ implies that 
		\begin{equation*}
			\begin{split}
				\int_{M}  Ric _{f-}^{H} \phi^{2} \lambda^{2} \mathrm{d}v_g 
				&\leq \epsilon\left\{10^{2 n}  R^{-1}vol(B_{3R})^{\frac{1}{n}} \right\}\int_{M}|\nabla(\phi \lambda)^{2}| \mathrm{d}v_g \\
				& \leq \epsilon\left\{10^{2 n}  R^{-1}vol(B_{3R})^{\frac{1}{n}} \right\}\bigg[\int_M \frac{1}{\epsilon}|\nabla \phi |^2\lambda^2+\epsilon \phi^2\lambda^2+ \epsilon_1|\nabla \lambda |^2\phi^2+\frac{1}{\epsilon_1}\phi^2\lambda^2 \mathrm{d} v_g\bigg],
			\end{split}
		\end{equation*}

		and we have( c.f. \cite{MR3692378} )	
				\begin{equation}\label{}
			\begin{split}
				&\int_{M} \phi^{2}\langle|d u|^{2 p-2} d u, \delta d\left(|d u|^{2 p-2} d u\right)\rangle\mathrm{d} v_g\\
				&\leq 2 \frac{p-2}{p-1} \int_{M} \phi|\nabla \phi||\nabla f||f|\mathrm{d} v_g+\left(\frac{p-2}{p-1}\right)^{2} \int_{M} \phi^{2}|\nabla f|^{2} \mathrm{d} v_g.
			\end{split}
		\end{equation}

		Thus combining  the above estimates  yields that
		\begin{equation}\label{qq}
			\begin{split}
				&\int_{M} \frac{1}{2} \phi^{2} \Delta \lambda^{2}\mathrm{d} v_g\\
				&+\epsilon\left\{10^{2 n} e^{2 a} Rvol(B_{3R})^{-\frac{1}{n}} \right\}\int_M\bigg[ \frac{1}{\epsilon_2}|\nabla \phi |^2\lambda^2\mathrm{d} v_g+ \epsilon_1|\nabla \lambda |^2\phi^2\mathrm{d} v_g+\left(\epsilon_2+ \frac{1}{\epsilon_1}\right)\phi^2 \lambda^2\mathrm{d} v_g\bigg]\\
				\geq&-2 \frac{p-2}{p-1} \int_{M} \phi|\nabla \phi||\nabla \lambda||\lambda|\mathrm{d} v_g-\left(\frac{p-2}{p-1}\right)^{2} \int_{M} \phi^{2}|\nabla \lambda|^{2}+\frac{n}{n-1} \int_{M} \phi^{2}|\nabla \lambda|^{2}\mathrm{d} v_g\\
				&+(n-1) H \int_{M} \lambda^{2} \phi^{2}\mathrm{d} v_g\\
				\geq&-2 \varepsilon_3 \int_{M} \phi^{2}|\nabla \lambda|^{2}\mathrm{d} v_g-\frac{1}{2 \epsilon_3}\left(\frac{p-2}{p-1}\right)^{2} \int_{M} \lambda^{2}|\nabla \phi|^{2}
				\mathrm{d} v_g-\left(\frac{p-2}{p-1}\right)^{2} \int_{M} \phi^{2}|\nabla \lambda|^{2}\mathrm{d} v_g\\
				&+\frac{n}{n-1} \int_{M} \phi^{2}|\nabla \lambda|^{2}\mathrm{d} v_g+(n-1) H \int_{M} \lambda^{2} \phi^{2}	  \mathrm{d}v_g, 
			\end{split}
		\end{equation}

		where the constant $p \geq 2$. Integrating by parts, we have 
		\begin{equation}\label{ww}
			\begin{split}
				\int_{M} \frac{1}{2} \phi^{2} \Delta \lambda^{2}\mathrm{d} v_g=-2 \int_{M} \phi|\lambda|\langle\nabla \phi, \nabla \lambda\rangle\mathrm{d} v_g \leq \delta\int_{M} |\phi|^{2}|\nabla \lambda|^{2}\mathrm{d} v_g+C_{\delta}\int_{M} |\nabla \phi|^{2}|\lambda|^{2}\mathrm{d} v_g.
			\end{split}
		\end{equation}
		
		Combining \eqref{ww} and \eqref{qq}, we have 
		\begin{equation}\label{}
			\begin{split}
				\left(\frac{n}{n-1}-\delta-2 \epsilon_3-\left(\frac{p-2}{p-1}\right)^{2}-\epsilon\epsilon_1\left\{10^{2 n}  R^{-1}vol(B_{3R})^{-\frac{1}{n}} \right\}\right) \int_{M} \phi^{2}|\nabla \lambda|^{2}\mathrm{d} v_g\\
				\leq\left(C_{\delta}+\epsilon\frac{1}{\epsilon_2}\left\{10^{2 n}  Rvol(B_{3R})^{-\frac{1}{n}} \right\}+\frac{1}{2 \epsilon_3}\left(\frac{p-2}{p-1}\right)^{2}\right) \int_{M} \lambda^{2}|\nabla \phi|^{2}\mathrm{d} v_g\\
				+	\left( \epsilon\left\{10^{2 n}  R^{-1}vol(B_{3R})^{-\frac{1}{n}} \right\}\left(\epsilon_2+ \frac{1}{\epsilon_1}\right)	- (n-1) H\right)  \int_{M} \lambda^{2} \phi^{2}	  \mathrm{d}v_g . 
			\end{split}
		\end{equation}
		Firstly,	we  choose sufficiently small $\epsilon_1,\epsilon_2 $,such that 	
		\begin{equation*}
			\begin{split}
				\epsilon\left\{10^{2 n} e^{2 a} R^{-1}vol(B_{3R})^{-\frac{1}{n}} \right\}\left(\epsilon_2+ \frac{1}{\epsilon_1}\right)	- (n-1) H \leq 0.
			\end{split}
		\end{equation*}
		
		Scondly, when $ R $ is sufficiently large,  we choose sufficiently small $\epsilon_3, \delta, $,such that 		
		$$ \frac{n}{n-1}-\delta-2 \epsilon_3-\left(\frac{p-2}{p-1}\right)^{2}-\epsilon\epsilon_1\left\{10^{2 n}  R^{-1}vol(B_{3R})^{-\frac{1}{n}} \right\}>0, $$
		
		then we can easily deduce that $|du|$ is a constant.  By \eqref{ev}, we see that $ \mathrm{Ric} (\omega^*,\omega^*) =0.$  	 We conclude that  $u$ is a constant.
	\end{proof}

	Motivated by \cite{wang2006harmonic}, we have 
	\begin{thm}
		let $ N^n $ be ompact  Riemannian manifold with nonpositive  sectional curvature.	Let $(M,g)$ be an $m(\geq 3)$-dimensional complete non-compact  metric measure space  with $ \lambda_1(\Delta)>0  $ and 
		\begin{equation*}
			\begin{split}
				Ric\geq -a\lambda_1(\Delta)+\delta,
			\end{split}
		\end{equation*} 
		for some $a\geq 0, \delta>0. $
		If u is  $p$-harmonic map with finite $L^{2p-2}$-energy from $ M $ to $ N $, then $ u $  is a constant.
	\end{thm}
	\begin{rem}
		This Theorem gerneralizes Theorem 1.1 in \cite{wang2006harmonic}.
	\end{rem}

	{\bf Proof}~~
	First we recall  Bochner formula  for $p$-harmonic map:
	for any smooth function $u$, we have (c.f. \cite[(3.13)]{MR3692378})
	\begin{align} \label{2.1.1}
		&\frac{1}{2}\Delta|du|^{2p-2} =|\nabla|du|^{p-2}du|^{2}-\langle |du|^{2p-2}du,\Delta\left( |du|^{p-2}du\right) \rangle \nonumber\\
		&+|du|^{2p-4}\langle du(Ric^{M}(e_{k}),du(e_{k})\rangle,
	\end{align}
	and the following  Kato inequality for p-harmonic map (cf.  \cite{MR3692378}):
	\begin{equation}\label{2.1.2}
		|\nabla\left(|du|^{p-2}  
		du\right) |^{2}\geq |\nabla|du|^{p-1}|^{2}.
	\end{equation}
	
	If we set $f=|du|^{p-1}$,  we have
	\begin{align}\label{2.1.4}
		\frac{1}{2}\Delta f^{2}+\bigg[-(1+\frac{1}{2ns})\lambda_1(M)+\delta\bigg]f^{2}
		\geq |\nabla f|^{2}-\langle|du|^{2p-2}du,\delta d(|du|^{2p-2}du)\rangle.
	\end{align}

	Now we choose a cut-off function in \eqref{f9} on $M$  , which satisfies that:
	\begin{equation}\label{e9}
		\begin{cases}
			0\leq\phi(x)\leq1 ,& x\in M \\ \phi(x)=1, & x\in B_{r}(x_{0}) \\ \phi(x)=0 ,&x\in M-B_{3r}(x_{0})\nonumber\\
			| \nabla \phi(x)|\leq \frac{C}{r}, & x \in M.
		\end{cases}
	\end{equation}
	Integrating on $M$ after multiplying both sides of (\ref{2.1.4}) by $\phi^{2}$, we get
	\begin{equation}\label{e2}
		\begin{split}
			&-\int_{M}\frac{1}{2}\phi^{2}\Delta f^{2}\mathrm{d} v_g\\
			\leq& \bigg[-a\lambda_1(M)+\delta\bigg]\int_{M}f^{2}\phi^{2}\mathrm{d} v_g -\int_{M}\phi^{2}|\nabla f|^{2}\mathrm{d} v_g\\
			&+\int_{M}\phi^{2}\langle |du|^{2p-2}du,\delta d(|du|^{2p-2}du)\rangle\mathrm{d} v_g.
		\end{split}
	\end{equation}

	Thus we  can cite the following inequality (c.f. \cite[(3.8)]{MR3692378} ),
	\begin{align}\label{2.1.8}
		&\int_{M}\phi^{2}\langle|du|^{2p-2}du,\delta d(|du|^{2p-2}du)\rangle\mathrm{d} v_g\\
		&\leq 2\frac{p-2}{p-1} \int_{M}\phi|\nabla\phi||\nabla f||f|
		\mathrm{d} v_g+(\frac{p-2}{p-1})^{2} \int_{M}\phi^{2}|\nabla f|^{2}\mathrm{d} v_g\\
		&\leq 2\varepsilon \int_{M}\phi^{2}|\nabla f|^{2}  \mathrm{d} v_g+\frac{1}{2\epsilon}(\frac{p-2}{p-1})^{2}\int_{M}f^{2}|\nabla \phi|^{2} \mathrm{d} v_g+(\frac{p-2}{p-1})^{2} \int_{M}\phi^{2}|\nabla f|^{2}\mathrm{d} v_g,
	\end{align}
	Hence ,we get 
	\begin{equation}\label{e2}
		\begin{split}
			&-\int_{M}\frac{1}{2}\phi^{2}\Delta f^{2}\mathrm{d} v_g\\
			\leq& \bigg[-a\lambda_1(M)+\delta\bigg]\int_{M}f^{2}\phi^{2}\mathrm{d} v_g \\
			&\frac{1}{2\epsilon}(\frac{p-2}{p-1})^{2}\int_{M}f^{2}|\nabla \phi|^{2}\mathrm{d} v_g +\left( (2\epsilon+\frac{p-2}{p-1})^{2}-1\right)  \int_{M}\phi^{2}|\nabla f|^{2}\mathrm{d} v_g,
		\end{split}
	\end{equation}
	Notice also that 
	\begin{equation}\label{o1}
		\begin{split}
			-\int_{M}\frac{1}{2}\phi^{2}\Delta f^{2}=2\int_M \phi f \nabla f\nabla \phi \leq \frac{1}{l}\int_M |\nabla f |^2 \phi^2+ l \int_M f |\nabla \phi |^2.
		\end{split}
	\end{equation}
	From the proof of Theorem 1.1 in \cite{1}, we know that:
	\begin{align}\label{874}
		\lambda_1(M)\int_{M}f^{2}\phi^{2}\mathrm{d} v_g\leq \int_{M}f^{2}|\nabla \phi|^{2}\mathrm{d} v_g+  \int_{M}\phi^{2}|\nabla f|^{2}\mathrm{d} v_g
		-\int_{M} \frac{1}{2}\phi^{2}\Delta f^{2}\mathrm{d} v_g.
	\end{align}
	
Combinging \eqref{874} and \eqref{e2} , we get 
	\begin{equation*}
		\begin{split}
			2 \int_{M^{n}} f \phi \nabla f \nabla \phi \mathrm{d} v_g
			=&\left(b_1+b_2\right) \int_{M^{n}} f \phi \nabla f\nabla \phi \mathrm{d} v_g\\
			\leq & 2b_1\bigg[\left\{\left( (2\epsilon+\frac{p-2}{p-1})^{2}-1\right)  \int_{M^{n}} \phi^{2}|\nabla f|^{2}\right.\left.\mathrm{d} v_g+\left(a \lambda_{1}-\delta\right) \int_{M^{n}} f^2 \phi^{2}\mathrm{d} v_g\right\} \\
			&+\frac{1}{2\epsilon}(\frac{p-2}{p-1})^{2}\int_{M}f^{2}|\nabla \phi|^{2}\mathrm{d} v_g\bigg]\\
			&+2b_2\left\{\frac{1}{l} \int_{M^{n}}\phi^{2}|\nabla f|^{2} \mathrm{d} v_g+l \int_{M^{n}} f^2|\nabla \phi|^{2}\mathrm{d} v_g\right\} \\
			=& 2b_1\left(a \lambda_{1}(M)-\delta\right) \int_{M^{n}} f^2 \phi^{2}\mathrm{d} v_g\\
			&+2\bigg[b_1\left( (2\epsilon+\frac{p-2}{p-1})^{2}-1\right)+\frac{b_2}{l}\bigg]  \int_{M^{n}} \phi^{2}|\nabla f|^{2} \mathrm{d} v_g\\
			&+2b_2 l\int_{M^{n}} f^2|\nabla \phi|^{2}+\frac{1}{2\epsilon}(\frac{p-2}{p-1})^{2}\int_{M}f^{2}|\nabla \phi|^{2}\mathrm{d} v_g,
		\end{split}
	\end{equation*}
	where $ b_1+b_2=1. $
	
	Plugging the above inequality into \eqref{874}, we get 
	
	\begin{equation}\label{}
		\begin{split}
			\lambda_1(M)\int_{M}f^{2}\phi^{2}\mathrm{d} v_g\leq& \int_{M}f^{2}|\nabla \phi|^{2}\mathrm{d} v_g+  \int_{M}\phi^{2}|\nabla f|^{2}\mathrm{d} v_g\\
			+& 2b_1\left(a \lambda_{1}(M)-\delta\right) \int_{M^{n}} f^2 \phi^{2}\mathrm{d} v_g\\
			&+2\bigg[b_1\left( (2\epsilon+\frac{p-2}{p-1})^{2}-1\right)+\frac{b_2}{l}\bigg]  \int_{M^{n}} \phi^{2}|\nabla f|^{2}\mathrm{d} v_g \\
			&+2b_2 l\int_{M^{n}} f^2|\nabla \phi|^{2}++\frac{1}{2\epsilon}(\frac{p-2}{p-1})^{2}\int_{M}f^{2}|\nabla \phi|^{2}\mathrm{d} v_g,
		\end{split}
	\end{equation}
	which is  written as 	
	\begin{equation}\label{}
		\begin{split}
			&\left( \lambda_1(M)-2ab_1 \lambda_{1}(M) +2b_1\delta\right) \int_{M}f^{2}\phi^{2}\mathrm{d} v_g\\
			\leq& \int_{M}f^{2}|\nabla \phi|^{2}\mathrm{d} v_g\\
			&+2\bigg[b_1\left( (2\epsilon+\frac{p-2}{p-1})^{2}-1\right)+\frac{b_2}{l}+1\bigg]  \int_{M^{n}} \phi^{2}|\nabla f|^{2} \mathrm{d} v_g\\
			&+2b_2 l\int_{M^{n}} f^2|\nabla \phi|^{2}\mathrm{d} v_g+\frac{1}{2\epsilon}(\frac{p-2}{p-1})^{2}\int_{M}f^{2}|\nabla \phi|^{2}\mathrm{d} v_g,
		\end{split}
	\end{equation}
	
	we can choose $ b_1  $ such that 
	\begin{equation*}
		\begin{split}
			\frac{-(1+\frac{1}{l}))}{2\epsilon+\frac{(p-2)^2}{(p-1)^2}-1-\frac{1}{l}}<b_1<\min\{1, \frac{1}{2a}\}
		\end{split}
	\end{equation*}
	
	Thus, we have $ f=0, $ thus, $ u $ is constant.

	\section{Liouville theorem in terms of BiRic curvature }

	\begin{thm}\label{52}
		Let $M$ be an $n(\geq 3)$-dimensional complete  stable  non-compact $ f $-minimal hypersurface  isometrically immersed in an $(n + 1)$-dimensional metric measure space $(\bar{M},g)$.  If  $\bar{M}$    has $\overline{BiRic}^{\delta}\geq 0$ with $ \delta \geq 1 $  and there exists at least a
		point $ x_0\in M $ such that $\overline{BiRic}^{\delta}(x_0)> 0,$  moreover, assume that $ \frac{1}{n-1}>\left( \frac{p-2}{p-1}\right)^2 $ then any  $p$-harmonic function with finite $L^{2p-2}$-energy on $M$  vanishes. If $ \min\{1, \frac{(p-1)^2}{(n-1)^2}\} >\left( \frac{p-2}{p-1}\right)^2 $, then any  $f$ wighted $p$-harmonic $ 1 $ forms with finite $L^{2p-2}$ energy on $M$  vanishes
	\end{thm}
	\begin{rem}
		This theorem generalize Theorem in \cite{zbMATH07153791} and Theorem in  \cite{zbMATH01208191}.
	\end{rem}
	\begin{proof}

		Since $M$ is $\delta$-stable,
		$$I(h)=\int_{M} \left( |\nabla \phi |^2-\delta(|A |^2+\overline{Ric_f}(n))\phi^2 \right)   \mathrm{d} v_g\geq 0.$$
		Take $ \phi=h|\omega|^{p-1} $, we derive 	
		\begin{equation}\label{w7}
			\begin{split}
				I(h)=-\int_{M} | h|^2\bigg(	|\omega|^{p-1}\Delta_f|\omega|^{p-1} +\delta(|A |^2+\overline{Ric_f}(n)) |\omega|^{2p-2} \bigg)\d v_g  +\int_{M}  |\nabla h|^2 |\omega|^{2p-2} \mathrm{d} v_g.
			\end{split}
		\end{equation}	
		However,		
		\begin{equation}\label{c1}
			\begin{split}
				&|\omega|^{p-1}\Delta_f|\omega|^{p-1}\\
				=	&|\nabla (|\omega|^{p-2}\omega)|^{2}-|\nabla| \omega|^{p-1}|^{2}+|\omega|^{2p-4}\sum_{i=1}^{n}\left \langle \omega(\mathrm{Ric}_f^M(e_i))  , \omega (e_i) \right \rangle \\
				&-\left \langle |\omega|^{p-2}\omega ,\Delta_f (|\omega|^{p-2}\omega) \right \rangle .
			\end{split}
		\end{equation}
		where	weighted Hodge laplacian
		$$\Delta_f=\delta_f \mathrm{d}+ \mathrm{d}\delta_f. $$

		Plugging \eqref{c1}into  \eqref{w7} to infer that 
		
		\begin{equation}\label{k1}
			\begin{split}
				&\int_{M}  |\nabla h|^2 |\omega|^{2p-2}\mathrm{d} v_g \\
				&\geq \int_{M} | h|^2\bigg(	|\omega|^{p-1}\Delta_f|\omega|^{p-1} +(|A |^2+\overline{Ric_f}(n)) |\omega|^{2p-2} \bigg) \d v_g\\
				&\geq \int_{M} | h|^2\bigg(F(\omega)- \left \langle |\omega|^{p-2}\omega ,\Delta_f (|\omega|^{p-2}\omega) \right \rangle+\Phi_2 +(|A |^2+\overline{Ric_f}(n)) |\omega|^{2p-2} \bigg)\d v_g.
			\end{split}
		\end{equation}
		Since $ Ric(X,X)\geq - \delta(\overline{Ric}(n,n)+|A |^2) $	
		\begin{equation}\label{k2}
			\begin{split}
				&\Phi_2+(|A |^2+\overline{Ric}(n)) |\omega|^{2p-2}\\
				&=	|\omega|^{2p-4} \sum R_{i j}\omega^{i} \omega_{j }
				+\delta(|A |^2+\overline{Ric}(n)) |\omega|^{2p-2}\\
				& \geq 0.
			\end{split}
		\end{equation}
		
		Combining \eqref{k1} and \eqref{k2}, we get 		
		\begin{equation}
			\begin{split}
				&\int_{M}	 h^2 	|\nabla (|\omega|^{p-2}\omega)|^{2}-|\nabla| \omega|^{p-1}|^{2} + \left \langle |\omega|^{p-2}\omega ,\Delta_f (|\omega|^{p-2}\omega) \right \rangle \d v_g \\
				\leq &\int_{M}  |\nabla h|^2 |\omega|^{2p-2} \mathrm{d} v_g.
			\end{split}
		\end{equation}

		Thus, by Kato's inequality for $ p $ harmonic function, 
		\begin{equation*}
			\begin{split}
				|\nabla (|du|^{p-2}du) |^2\geq \frac{n}{n-1}|\nabla (|du|^{p-1}) |^2.
			\end{split}
		\end{equation*}
		we have 
		\begin{equation}
			\begin{split}
				\int_{M}	\frac{1}{n-1}|\nabla| \omega|^{p-1}|^{2} \mathrm{d} v_g 
				\leq& \int_{M}  |\nabla h|^2 |\omega|^{2p-2} \mathrm{d} v_g-\left \langle |\omega|^{p-2}\omega ,\Delta (|\omega|^{p-2}\omega) \right \rangle h^2\mathrm{d} v_g\\
				\leq& \int_{M}  |\nabla h|^2 |\omega|^{2p-2} \mathrm{d} v_g+2\frac{p-2}{p-1}\int_{M}h|\nabla h||\nabla |du|^{p-1}||du|^{p-1}\mathrm{d} v_g\\
				&+\left( \frac{p-2}{p-1}\right)^2\int_{M}h^2|\nabla |du|^{p-1} |^2 \mathrm{d} v_g.
			\end{split}
		\end{equation}
		Hence, we deduce that 
		
		\begin{equation}
			\begin{split}
				&\left(\frac{1}{n-1}- \left( \frac{p-2}{p-1}\right)^2	-\epsilon\right)	\int_{M} |\nabla| \omega|^{p-1}|^{2}\mathrm{d} v_g\\ 
				\leq& \bigg[1+\frac{1}{4\epsilon}4 \left( \frac{p-2}{p-1}\right)^2\bigg]\int_{M}  |\nabla h|^2 |\omega|^{2p-2} \mathrm{d} v_g\\
			\end{split}
		\end{equation}
		By \eqref{c1}, we infer that  $ |\omega| $ is  constatn. We may assume  $ |\omega|\neq 0,$ otherwise the theorem has been proved.

		By \eqref{k1}, we have 
		
		\begin{equation*}
			\begin{split}
				Ric(\omega^*,\omega^*)=0
			\end{split}
		\end{equation*}
		
		By \eqref{k1} and \eqref{k2},we have 
		
		\begin{equation*}
			\begin{split}
				Ric(\omega^*,\omega^*)	+\delta(|A |^2+\overline{Ric}(n)) |\omega|^{2p-2}=0.
			\end{split}
		\end{equation*}
		However, by \cite{zbMATH07153791}, we know that 
		\begin{equation*}
			\begin{split}
				Ric(X,X)\geq \overline{BiRic}^{\delta}(X,N)-\delta((|A |^2+\overline{Ric_f}(n)))
			\end{split}
		\end{equation*}
		Thus we have $ \overline{BiRic}=0 $, this is a contradiction to our assumption in the theorem.
		The last statement follows from  the kato inequality for $ p $ harmonic 1 form.
	\end{proof}

	\bigskip

	Follow the proof of  Cheng \cite{zbMATH01208191} and 
	Li and Zou \cite{zoulibiricci}, for constant mean curvature submanifold, we  can get
	
	\begin{thm}\label{8t}
		Let M be an n-dimensional complete and non
		compact orientable hypersurface with constant mean curvature  $  H $ in a Riemannian
		manifold with biRicci curvature satisfying along M 
		\begin{equation*}
			\begin{split}
				\overline{BiRic}^\delta\geq \frac{(n-5)n^2}{4}H^2.
			\end{split}
		\end{equation*}If $ M $ is strongly  $\delta$-stable and $ p\geq 2 , \frac{1}{n-1}>\left( \frac{p-2}{p-1}\right)^2$, then there are
		no nontrivial $ L^{2p-2} $ $ p $-harmonic function  on $ M $. Morover, If $ 2\leq p<2+\frac{1}{\sqrt{n-1}} $, then there are
		no nontrivial $ L^p $ $ p $-harmonic $ 1 $-forms on $ M $. Generally, if $ \min\{1, \frac{(p-1)^2}{(n-1)^2}\}>(p-2)^2 $, there are no nontrivial $ L^p $ $ p $-harmonic $ 1 $-forms on $ M $
	\end{thm}
	\begin{rem}
		When p=2, it is just Theorem 1 in Cheng \cite{zbMATH01208191}.
	\end{rem}
	\begin{rem}
		The proof depends on the good bound in Kato's inequality, thus the theorem can not be applied to $ p $ harmonic map or $  p $ harmonic $ k $ form $ (k\geq 2) $.
	\end{rem}
	\begin{proof}  We slightly modify the proof of Theorem  \ref{52}. Let $ \omega=du $, 		
		\begin{equation}\label{y2}
			\begin{split}
				|\omega|^{p-1}\Delta|\omega|^{p-1}=F(\omega)+\Phi_2- \left \langle |\omega|^{p-2}\omega ,\Delta(|\omega|^{p-2}\omega) \right \rangle .
			\end{split}
		\end{equation}
		wherer $ \Delta =d\delta+\delta d,$ $ \Phi_2=	|\omega|^{2p-4} \sum R_{i j}\omega^{i} \omega_{j }. $ 
		
		Since $ M $ is $ \delta $ stable, 
		\begin{equation*}
			\begin{split}
				&\int_{M}  |\nabla h|^2 |\omega|^{2p-2}\mathrm{d} v_g \\
				&\geq \int_{M} | h|^2\bigg(	|\omega|^{p-1}\Delta|\omega|^{p-1} +\delta(|A |^2+\overline{Ric}(n)) |\omega|^{2p-2} \bigg)\mathrm{d} v_g \\
				&\geq \int_{M} | h|^2\bigg(F(\omega)+ \left \langle |\omega|^{p-2}\omega ,\Delta (|\omega|^{p-2}\omega) \right \rangle+\Phi_2 +\delta(|A |^2+\overline{Ric}(n)) |\omega|^{2p-2} \bigg)\mathrm{d} v_g .
			\end{split}
		\end{equation*}
		where $ F(\omega)=|\nabla (|\omega|^{p-2}\omega)|^{2}-|\nabla| \omega|^{p-1}|^{2}. $ 
		
		Notice also that 
		\begin{equation}\label{3.10}
			\begin{split}
				&\Phi_2+\delta(|A |^2+\overline{Ric}(n)) |\omega|^{2p-2}\\
				=&	|\omega|^{2p-4} \sum R_{i j} \omega^{i} \omega_{j }
				+\delta(|A |^2+\overline{Ric}(n)) |\omega|^{2p-2}\\
				=& 	|\omega|^{2p-4} \sum \bigg( \overline{Ric}_{ij}-\bar{R}(n,e_j,n,e_i) + nH A_{ij} -\sum_{k}A_{ik}A_{jk} \bigg) \omega^{i} \omega_{j }
				\\
				&+\delta(|A |^2+\overline{Ric}(n)) |\omega|^{2p-2}\\
			\end{split}
		\end{equation}
		
		By \cite[Lemma 2.1]{zbMATH01208191}, we have 
		\begin{equation*}
			\begin{split}
				&\left( nH A_{ij} -\sum_{k}A_{ik}A_{jk}\right)\omega^i\omega_j+|A|^2|\omega|^2 \\
				=&n H A\left(\omega^{*}, \omega^{*}\right)-\left\langle A \omega^{*}, A \omega^{*}\right\rangle+|A|^{2}|\omega|^{2}\geq -\frac{(n-5)n^2 H^2}{4} |\omega|^2.
			\end{split}
		\end{equation*}
		Since $ \mathrm{BiRic}^\delta \geq \frac{(n-5)n^2 H^2}{4}$, we have 
		\begin{equation*}
			\begin{split}
				\Ric(X,X)+\delta \Ric(n)|X |^2-K(X,n)|X |^2\geq  \frac{(n-5)n^2 H^2}{4} |X |^2.
			\end{split}
		\end{equation*}
		where $ X=\omega^{*}. $
		\begin{equation}\label{y1}
			\begin{split}
				&\int_{M}  |\nabla h|^2 |\omega|^{2p-2} \mathrm{d} v_g \\
				&\geq \int_{M} | h|^2\bigg(	|\omega|^{p-1}\Delta|\omega|^{p-1} +\delta(|A |^2+\overline{Ric}(n)) |\omega|^{2p-2} \bigg)\mathrm{d} v_g \\
				&\geq \int_{M} | h|^2\bigg(F(\omega)-\left \langle |\omega|^{p-2}\omega ,\Delta (|\omega|^{p-2}\omega) \right \rangle)  + BiRic (X,n)-\frac{(n-5)n^2 H^2}{4} |\omega|^2\bigg)\mathrm{d} v_g .
			\end{split}
		\end{equation}
		By Kato inequality for $ p $ harmonic function, we have 
		\begin{equation}
			\begin{split}
				&	\frac{1}{n-1}\int_{M} h^2	|\nabla| \omega|^{p-1}|^{2} \mathrm{d} v_g \\ 
				\leq& \int_{M}  |\nabla h|^2 |\omega|^{2p-2} \mathrm{d} v_g-\left \langle |\omega|^{p-2}\omega ,\Delta (|\omega|^{p-2}\omega) \right \rangle h^2\mathrm{d} v_g \\
				\leq& \int_{M}  |\nabla h|^2 |\omega|^{2p-2} \mathrm{d} v_g+2\frac{p-2}{p-1}\int_{M}h|\nabla h||\nabla |du|^{p-1}||du|^{p-1}\mathrm{d} v_g \\
				&+\left( \frac{p-2}{p-1}\right)^2\int_{M}h^2|\nabla |du|^{p-1} |^2 \mathrm{d} v_g .
			\end{split}
		\end{equation}
		Thus, we have 	
		\begin{equation}
			\begin{split}
				\left( \frac{1}{n-1}-\left( \frac{p-2}{p-1}\right)^2-\epsilon\right)	\int_{M} |\nabla| \omega|^{p-1}|^{2}  \mathrm{d} v_g  
				\leq& C(p)\int_{M}  |\nabla h|^2 |\omega|^{2p-2} \mathrm{d} v_g.\\
			\end{split}
		\end{equation}
		By kato's inequyality for  $ p $ harmonic 1 form $ (p \geq 2) $, from \eqref{y1}, we see that 
		\begin{equation}
			\begin{split}
				\left( \frac{1}{(n-1)(p-1)^2}-\left( \frac{p-2}{p-1}\right)^2-\epsilon\right)	\int_{M} |\nabla| \omega|^{p-1}|^{2}  \mathrm{d} v_g  
				\leq& C(p)\int_{M}  |\nabla h|^2 |\omega|^{2p-2} \mathrm{d} v_g.\\
			\end{split}
		\end{equation}
		By \eqref{y2}, we have 		
		\begin{equation*}
			\begin{split}
				Ric(\omega^*,\omega^*)=0
			\end{split}
		\end{equation*}		
		By  \eqref{y1},we have 		
		\begin{equation*}
			\begin{split}
				Ric(\omega^*,\omega^*)	+\delta(|A |^2+\overline{Ric}(n)) |\omega|^{2p-2}=0.
			\end{split}
		\end{equation*}
		However, by  Gauss equation and the definition of $ \mathrm{BiRic}^\delta $, we know that (c.f. \cite{zbMATH07153791} and \cite[Theorem 1]{zbMATH01208191}) 
		\begin{equation*}
			\begin{split}
				\Ric(X,X)=&\overline{\mathrm{Ric}}(X,X)+\delta\overline{\mathrm{Ric}}(n,n)|X |^2-K(X,n)|X |^2\\
				&-\delta\left( \overline{\mathrm{Ric}}(n,n)+|A |^2\right)|X |^2+ nH \langle AX,X\rangle-|AX |^2+a|A |^2|X |^2 \\
				\geq& \overline{\mathrm{BiRic}}^{\delta}(\frac{X}{|X|},N)|X |^2-\frac{(n-5)n^2 H^2}{4} |X|^2-\delta((|A |^2+\overline{\Ric}(n)))|X |^2
			\end{split}
		\end{equation*}

		As in  Cheng \cite[Theorem 1]{zbMATH01208191}, we can conclude that 
		\begin{equation*}
			\begin{split}
				\mathrm{Ric}_M \geq 0.
			\end{split}
		\end{equation*}
		It is well known that  $ M $ has infinite volume (c.f.  Yau\cite{yaust}). Since $ |\omega|  $ is a  positive constant, this is  a contradiction to the fact that $\omega$ has finite $ L^{2p-2} $ energy.
		
		The last statement follows from the Kato inequality in Lemma \ref{kato} and 
		\begin{equation}
			\begin{split}
				\left( \min\{1, \frac{(p-1)^2}{(n-1)^2}\}-\left( \frac{p-2}{p-1}\right)^2-\epsilon\right)	\int_{M} |\nabla| \omega|^{p-1}|^{2}  \mathrm{d} v_g  
				\leq& C(p)\int_{M}  |\nabla \eta|^2 |\omega|^{2p-2} \mathrm{d} v_g.\\
			\end{split}
		\end{equation}
		
	\end{proof}

	\section{Liouville type theorem for  harmonic $ q $ form}

	Inspired by the method in \cite{zbMATH01208191}, we get 
	\begin{thm}
		Let $M$ be an $n(2< n\leq 4)$-dimensional complete 2-stable  noncompact minimal hypersurface  isometrically immersed in an $(n + 1)$-dimensional  a complete noncompact metric measure space  $(\bar{M},g), $ which satisifies and $ \overline{\mathrm{BiRic}}\geq 0. $  Moreover, assume that $ \overline{R}_{k j i h} T^{kj}T^{ih}\leq 0$ for any antisymmetric tensor $ T $ ,  then any   harmonic $ 2 $-form with finite $L^2$-energy on $M$  vanishes.
	\end{thm}
	\begin{rem}
		The theorem can applied to the case where $ \bar{M}=\mathbb{R}^{n+1} $ or  $ \mathbb{R}^{n+1} $, thus we generalize Theorem 6.1 in \cite{zbMATH01208191}.
	\end{rem}
	\begin{proof}
		
		For  any differential $q$ form  $ \omega $,  it is well known that 
		\[	
		\Delta \omega=\Delta \omega_{i_{1} \cdots i_{q}} =\nabla^{r} \nabla_{r} \omega_{i_{1} \cdots l_{q}}-\sum_{s=1}^{q} R_{i_{s}}^{r} \omega_{i_{1} \cdots r \cdots i_{q}}+\sum_{t<s}^{ q} R_{i_ti_{s}}^{v u} \omega_{i_{1} \cdots v \cdots u \cdots i_{q}} \\
		=0 .
		\]
		Putting $\|\omega\|^{2}=\sum \omega_{i_{1} \cdots i_{q}} \omega^{l_{1} \cdots i_{q}}$ and $\|\nabla \omega\|^{2}=\sum \nabla_{r} \omega_{i_{1} \cdots i_{q}} \nabla^{r} \omega^{i_{1} \cdots i_{q}}$, we obtain (c.f. \cite{zbMATH01208191})
		\begin{equation}
			\begin{split}
				\frac{1}{2} \Delta\|\omega\|^{2} &=\|\nabla \omega\|^{2}+\sum \omega_{i_{1} \cdots l_{q}} \nabla^{r} \nabla_{r} \omega^{i_{1} \cdots i_{q}} \\
				&=\|\nabla \omega\|^{2}+\sum R_{i_{s}}^{r} \omega_{i_{1} \cdots r \cdots i_{q}} \omega^{i_{1} \cdots i_{q}}-\sum_{t<s}^{1 \cdots q} R_{{i_t} i_{s}}^{v u} \omega_{i_{1} \cdots v \cdots u \cdots i_{p}} \omega^{i_{1} \cdots i_{q}} \\
				&	=\|\nabla \omega\|^{2}+q \sum R_{i_{1}}^{r} \omega_{ri_{2} \cdots i_{q}} \omega^{i_{1} \cdots i_{q}}-\sum_{t<s}^{1 \cdots p} R_{i_ti_s}^{v u} \omega_{i_{1} \cdots v \cdots u \cdots i_{q}} \omega^{i_{1} \cdots i_{q}} \\
				&	=\|\nabla \omega\|^{2}+q \sum R_{i j} \omega_{i_{2} \cdots i_{q}}^{i} \omega_{j i_{2} \cdots i_{q}}-\frac{q(q-1)}{2} \sum R_{t j i h} \omega_{i_{3} \cdots i_{q}}^{tj} \omega^{i h{i_{3} \cdots i_{q}}} +\Delta \omega.
			\end{split}
		\end{equation}
		%

		On the other hand, we have
		\[
		\frac{1}{2} \Delta\|\omega\|^{2} =\|\omega\| \Delta\|\omega\|+\|\nabla\| \omega\|\|^{2} \\
		=\|\omega\| \Delta\|\omega\|+\|\nabla \omega\|^{2}-F(\omega),
		\]
		where
		\[
		F(\omega)=\|\nabla \omega\|^{2}-\|\nabla\| \omega\|\|^{2}
		\]

		We take cut off function $ h(x) $ as  in \eqref{f9}	,
		\begin{equation}\label{cf1}
			\begin{cases}
				0\leq h(x)\leq1 ,& x\in M 
				\\ h(x)=1, & x\in B_{r}(x_{0}) \\ 
				h(x)=0 ,&x\in M-B_{2r}(x_{0})\\
				| \nabla h|\leq \frac{C}{r}, & x \in M.
			\end{cases}
		\end{equation}
		
		Since $ M $ is $  2  $-stable, we have 
		\begin{equation}
			\begin{split}
				&\int_{M}  |\nabla h|^2 |\omega|^{2} \mathrm{d} v_g \\
				&\geq \int_{M} | h|^2\bigg(	|\omega|\Delta|\omega| +2(|A |^2+\overline{Ric}(n)) |\omega|^{2} \bigg)\mathrm{d} v_g\\
				&\geq \int_{M} | h|^2\bigg(F(\omega)+\Phi_2 +2(|A |^2+\overline{Ric}(n)) |\omega|^{2} \bigg)\mathrm{d} v_g.
			\end{split}
		\end{equation}
		where 
		\begin{equation}\label{3.9}
			\begin{split}
				\Phi_2=&2 \sum R_{i j} \omega_{i_{2} }^{i} \omega_{j i_{2} }
				- \sum R_{k j i h} \omega_{}^{k{ }{j}} \omega^{i h{}} \\
				=&	2 \sum   R_{i j}\omega_{i_{2} }^{i} \omega_{j i_{2} }
				-\sum R_{k j i h} \omega_{}^{k{ }{j}} \omega^{i h{}}.
			\end{split}
		\end{equation}
		Notice  that 
		\begin{equation}\label{385}
			\begin{split}
				&\Phi_2+2(|A |^2+\overline{Ric}(n)) |\omega|^{2}\\
				&=	2 \sum R_{i j} \omega_{i_{2} }^{i} \omega_{j i_{2} }
				- \sum R_{k j i h} \omega_{}^{k{ }{j}} \omega^{i h{}} \\
				&+2(|A |^2+\overline{Ric}(n)) |\omega|^{2p-2}\\
				&= 	2|\omega|^{2p-4} \sum \bigg(\left( \overline{Ric}\right)_{ij}-\bar{R}(n,e_j,n,e_i) + nH A_{ij} -\sum_{k}A_{ik}A_{jk}  \bigg) \omega_{i_{2} }^{i} \omega_{j i_{2} }
				\\
				&- \sum \big(\overline{R}_{k j i h}+A_{ki}A_{jh}-A_{kh}A_{ji}\big) \omega_{}^{k{ }{j}} \omega^{i h{}}
				+q(|A |^2+\overline{Ric}(n)) |\omega|^{2}.
			\end{split}
		\end{equation}
		
		%
		Since $ \mathrm{BiRic}_f \geq 0$, we have 
		\begin{equation*}
			\begin{split}
				Ric_f(X,X)+Ric_f(n)-K(X,n)\geq 0,
			\end{split}
		\end{equation*}
		where $ X=\frac{\omega^{*}}{|\omega^{*}|}. $

		
		By \cite[Lemma 6.1]{zbMATH01208191}, when $ 2\leq n \leq 4,  $ we have 		
		\begin{equation}\label{}
			\begin{split}
				2&\bigg(nH A_{ij} -\sum_{k}A_{ik}A_{jk} \bigg) \omega_{i_{2} }^{i} \omega_{j i_{2} }\\
				-&(A_{ki}A_{jh}-A_{kh}A_{ji}\big) \omega_{}^{k{ }{j}} \omega^{i h{}}+|A |^2|\omega|^2 \geq 0.
			\end{split}
		\end{equation}
		
		Moreover by our assumption, we have 
		\begin{equation*}
			\begin{split}
				-\overline{R}_{k j i h}\omega_{}^{k{ }{j}} \omega^{i h{}}\geq 0.
			\end{split}
		\end{equation*}
		Plugging all the above estimates into \eqref{385}, we have 		
		\begin{equation*}
			\begin{split}
				\Phi_2+(|A |^2+\overline{Ric}(n)) |\omega|^{2p-2}\geq 0.
			\end{split}
		\end{equation*}
		By Kato inequality for harmonic $ 2  $ form (c.f.\cite{calderbank2000refined} ), we have 
		\begin{equation}
			\begin{split}
				&\frac{1}{n-2}\int_{M} h^2	|\nabla| \omega|^2|^{2}\mathrm{d} v_g 
				\leq \int_{M}  |\nabla h|^2 |\omega|^{2} \mathrm{d} v_g
			\end{split}
		\end{equation}
		
		Thus, we have 	
		\begin{equation}
			\begin{split}
				\int_{M} h^2 |\nabla| \omega||^{2} \mathrm{d} v_g 
				\leq& (n-2)\int_{M}  |\nabla h|^2 |\omega|^{2} \mathrm{d} v_g\\
			\end{split}
		\end{equation}
		It is standard to infer that  $ \omega $ is parallel. One can refer to the argument in in \cite[page 9]{zbMATH01208191}.
	\end{proof}
	
	\begin{thm}
		Let $M$ be an $n(\geq 3)$-dimensional complete $ q $-stable  non-compact submanifold isometrically immersed in an $(n + k)$-dimensional meetric measure space $(\bar{M},g)$  If  $ \bar{M}$  has  pure curvature tensor  and $ \overline{\mathrm{BiRic}}\geq 0. $  Moreover, assume that $ \overline{R}_{k j i h} T^{kj}T^{ih}\leq 0$ for any antisymmetric tensor $ T $ , M has flat normal bundle. If one of the following conditons holds,

		(1)$ \frac{n}{3} \leq q, $
		
		(2) $ |A |^2 \leq \frac{n^2 |H |^2}{C_{n,q}}, $	
		
		then any harmonic $ q $-form with finite $L^{2}$-energy on $M$  must be parallel.
	\end{thm}
	\begin{rem}
		When $ \bar{M}=\mathbb{S}^{n+1} $ or $ \mathbb{R}^{n+1} $, the condition for $ \bar{M} $ holds. This theorem doesn't hold for $ p (\neq 2)$ harmonic $ q $ form because of the lack of good Kato's inequality.
	\end{rem}
	\begin{proof}
		By Guass equation, we have 
		\begin{equation}
			\begin{split}
				&	q\sum  R_{i j}\omega_{i_{2} \cdots i_{q}}^{i} \omega_{j i_{2} \cdots i_{q}}
				-\frac{q(q-1)}{2} \sum R_{k j i h} \omega_{i_{3} \cdots i_{q}}^{k{ }{j}} \omega^{i h{i_{3} \cdots i_{q}}} \\
				&+(|A |^2+\overline{Ric}(n)) |\omega|^{2p-2}\\
				&= 	q \sum \bigg( \overline{R}_{ikjk} + nH A_{ij} -\sum_{k}A_{ik}A_{jk} \bigg) \omega_{i_{2} \cdots i_{q}}^{i} \omega_{j i_{2} \cdots i_{q}}
				\\
				&-\frac{q(q-1)}{2} \sum \big(\overline{R}_{k j i h}+A_{ki}A_{jh}-A_{kh}A_{ji}\big) \omega_{i_{3} \cdots i_{q}}^{k{ }{j}} \omega^{i h{i_{3} \cdots i_{q}}}
				+q(|A |^2+\overline{Ric}(n)) |\omega|^{2p-2}.
			\end{split}
		\end{equation}
		Since $ \mathrm{BiRic}_f \geq 0$, we have 
		\begin{equation*}
			\begin{split}
				Ric(X,X)+Ric(n)-K(X,n)\geq 0.
			\end{split}
		\end{equation*}
		where $ X=\frac{\omega^{*}}{|\omega^{*}|}. $ Thus, we have 
		\begin{equation}\label{pp2}
			\begin{split}
				q \sum  \overline{R}_{ikjk}  \omega_{i_{2} \cdots i_{q}}^{i} \omega_{j i_{2} \cdots i_{q}}+q\overline{Ric}(n) |\omega|^{2}\geq 0.
			\end{split}
		\end{equation}
		
		By\cite{zbMATH06505276}, we have
		\begin{equation}
			\begin{split}
				&	q \sum \bigg( nH A_{ij} -\sum_{k}A_{ik}A_{jk}  \bigg) \omega_{i_{2} \cdots i_{q}}^{i} \omega_{j i_{2} \cdots i_{q}}
				\\
				&-|\omega|^{2p-4}\frac{q(q-1)}{2} \sum \big(A_{ki}A_{jh}-A_{kh}A_{ji}\big) \omega_{i_{3} \cdots i_{q}}^{k{ }{j}} \omega^{i h{i_{3} \cdots i_{q}}} \\
				&\geq  \left( \frac{1}{2}(n^2|H |^2-C_{n,q}|A |^2)\right).	
			\end{split}
		\end{equation}
		Thus, we have 			
		\begin{equation}\label{pp1}
			\begin{split}
				&	q \sum \bigg( nH A_{ij} -\sum_{k}A_{ik}A_{jk}  \bigg) \omega_{i_{2} \cdots i_{q}}^{i} \omega_{j i_{2} \cdots i_{q}}\\
				& 	-\frac{q(q-1)}{2} \sum \big(A_{ki}A_{jh}-A_{kh}A_{ji}\big) \omega_{i_{3} \cdots i_{q}}^{k{ }{j}} \omega^{i h{i_{3} \cdots i_{q}}} +	q|A |^2\\
				& \geq  \left( \frac{1}{2}(n^2|H |^2-C_{n,q}|A |^2) +	q|A |^2\right) \geq 0 .	
			\end{split}
		\end{equation}
		where in the last inequality, we used the conditions in the theorem.

		By \eqref{pp1} and \eqref{pp2}, we have 	
		\begin{equation}
			\begin{split}
				&\int_{M}  |\nabla h|^2 |\omega|^{2} \mathrm{d} v_g \\
				\geq& \int_{M} | h|^2\bigg(	|\omega|\Delta|\omega| +(|A |^2+\overline{Ric}(n)) |\omega|^{2} \bigg)\mathrm{d} v_g\\
				\geq& \int_{M} | h|^2\bigg(F(\omega)+ \Phi_2 +q(|A |^2+\overline{Ric}(n)) |\omega|^{2} \bigg)\mathrm{d} v_g\\
				\geq &\int_{M} | h|^2F(\omega)\mathrm{d} v_g .
			\end{split}
		\end{equation}
		
		By Kato inequality for harmonic $ q  $ form (c.f.\cite{calderbank2000refined} ), we have 
		\begin{equation}
			\begin{split}
				&K_q\int_{M} h^2	|\nabla| \omega||^{2} \mathrm{d} v_g 
				\leq \int_{M}  |\nabla h|^2 |\omega|^{2} \mathrm{d} v_g
			\end{split}
		\end{equation}
		where 
		\begin{equation*}
			\begin{split}
				K_q= \begin{cases}
					\frac{1}{n-q}&  if 2 \leq q \leq \frac{n}{2} \\
					\frac{1}{q}	&  if \frac{n}{2} \leq q \leq n-2,
				\end{cases}
			\end{split}
		\end{equation*}
		Thus, we have 	
		\begin{equation}
			\begin{split}
				K_q	\int_{M} h^2 |\nabla| \omega||^{2}   \mathrm{d} v_g
				\leq& \int_{M}  |\nabla h|^2 |\omega|^{2} \mathrm{d} v_g\\
			\end{split}
		\end{equation}
	Now it is standard to infer that  $ \omega $ is parallel.
	\end{proof}

		\bibliographystyle{plain}
		\bibliography{VT-libs-2022}

\begin{thebibliography}{10}

\bibitem{p-harmonicmap}
Paul Baird and Sigmundur Gudmundsson.
\newblock p-harmonic maps and minimal submanifolds.
\newblock {\em Mathematische Annalen}, 294(1):611--624, 1992.

\bibitem{Baird1992}
Paul Baird and Sigmundur Gudmundsson.
\newblock {$p$}-harmonic maps and minimal submanifolds.
\newblock {\em Math. Ann.}, 294(4):611--624, 1992.

\bibitem{zbMATH04182289}
Pierre B{\'e}rard.
\newblock A note on {Bochner} type theorems for complete manifolds.
\newblock {\em Manuscr. Math.}, 69(3):261--266, 1990.

\bibitem{calderbank2000refined}
David~MJ Calderbank, Paul Gauduchon, and Marc Herzlich.
\newblock Refined kato inequalities and conformal weights in riemannian
  geometry.
\newblock {\em Journal of Functional Analysis}, 173(1):214--255, 2000.

\bibitem{MR3692378}
Xiangzhi Cao.
\newblock Liouville type theorem about {$p$}-harmonic function and
  {$p$}-harmonic map with finite {$L^q$}-energy.
\newblock {\em Chinese Ann. Math. Ser. B}, 38(5):1071--1076, 2017.

\bibitem{zbMATH07361892}
Xiaoli Chao, Aiying Hui, and Miaomiao Bai.
\newblock Vanishing theorems for {{\(p\)}}-harmonic {{\(\ell\)}}-forms on
  {Riemannian} manifolds with a weighted {Poincar{\'e}} inequality.
\newblock {\em Differ. Geom. Appl.}, 76:13, 2021.
\newblock Id/No 101741.

\bibitem{cheng20002}
Xu~Cheng.
\newblock {$L^2$} harmonic forms and stability of hypersurfaces with constant
  mean curvature.
\newblock {\em Boletim da Sociedade Brasileira de
  Matem{\'a}tica-Bulletin/Brazilian Mathematical Society}, 31(2):225--239,
  2000.

\bibitem{cheng2006constant}
Xu~Cheng.
\newblock On constant mean curvature hypersurfaces with finite index.
\newblock {\em Archiv der Mathematik}, 86(4):365--374, 2006.

\bibitem{zbMATH02114455}
Xianzhe Dai and Guofang Wei.
\newblock A heat kernel lower bound for integral {Ricci} curvature.
\newblock {\em Mich. Math. J.}, 52(1):61--69, 2004.

\bibitem{zbMATH06668693}
Nguyen~Thac Dung.
\newblock {{\(p\)}}-harmonic {{\(\ell\)}}-forms on {Riemannian} manifolds with
  a weighted {Poincar{\'e}} inequality.
\newblock {\em Nonlinear Anal., Theory Methods Appl., Ser. A, Theory Methods},
  150:138--150, 2017.

\bibitem{zbMATH07153791}
Nguyen~Thac Dung, Nguyen~Van Duc, and Juncheol Pyo.
\newblock Harmonic 1-forms on immersed hypersurfaces in a {Riemannian} manifold
  with weighted bi-{Ricci} curvature bounded from below.
\newblock {\em J. Math. Anal. Appl.}, 484(1):29, 2020.
\newblock Id/No 123693.

\bibitem{dung2017p}
Nguyen~Thac Dung and Keomkyo Seo.
\newblock p-harmonic functions and connectedness at infinity of complete
  submanifolds in a riemannian manifold.
\newblock {\em Annali di Matematica Pura ed Applicata (1923-)},
  196(4):1489--1511, 2017.

\bibitem{MR3849353}
Nguyen~Thac Dung and Pham~Trong Tien.
\newblock Vanishing properties of {$p$}-harmonic {$\ell$}-forms on {R}iemannian
  manifolds.
\newblock {\em J. Korean Math. Soc.}, 55(5):1103--1129, 2018.

\bibitem{1}
Wenzhen Gan and Peng Zhu.
\newblock ${L}^2$ harmonic 1-forms on minimal submanifolds in spheres.
\newblock {\em Results in Mathematics}, 65(3-4):483--490, 2014.

\bibitem{zbMATH06859577}
Yingbo Han.
\newblock {{\(p\)}}-harmonic {{\(l\)}}-forms on complete noncompact
  submanifolds in sphere with flat normal bundle.
\newblock {\em Bull. Braz. Math. Soc. (N.S.)}, 49(1):107--122, 2018.

\bibitem{zbMATH07040685}
Yingbo Han.
\newblock Vanishing theorem for {{\(p\)}}-harmonic 1-forms on complete
  submanifolds in spheres.
\newblock {\em Bull. Iran. Math. Soc.}, 44(3):659--671, 2018.

\bibitem{zbMATH06451355}
Yingbo Han and Shuxiang Feng.
\newblock A {Liouville} type theorem for {{\(p\)}}-harmonic functions on
  minimal submanifolds in {{\(\mathbb R^{n+m}\)}}.
\newblock {\em Mat. Vesn.}, 65(4):494--498, 2013.

\bibitem{2}
Yingbo Han and Hong Pan.
\newblock ${L}^p$ $p$-harmonic 1-forms on submanifolds in a {H}adamard
  manifold.
\newblock {\em Journal of Geometry and Physics}, 107:79--91, 2016.

\bibitem{zbMATH06827101}
Yingbo Han, Qianyu Zhang, and Mingheng Liang.
\newblock {{\(L^{p}\)}} {{\(p\)}}-harmonic 1-forms on locally conformally flat
  {Riemannian} manifolds.
\newblock {\em Kodai Math. J.}, 40(3):518--536, 2017.

\bibitem{Kawai1999}
Shigeo Kawai.
\newblock p-{H}armonic {M}aps and {C}onvex {F}unctions.
\newblock {\em Geometriae Dedicata}, 74(3):261--265, 1999.

\bibitem{zbMATH01208191}
Haizhong Li.
\newblock {{\(L^2\)}} harmonic forms on a complete stable hypersurfaces with
  constant mean curvature.
\newblock {\em Kodai Math. J.}, 21(1):1--9, 1998.

\bibitem{Li2001}
Jin~Tang Li.
\newblock {$P$}-harmonic maps for submanifolds with positive {R}icci curvature.
\newblock {\em Xiamen Daxue Xuebao Ziran Kexue Ban}, 40(6):1191--1195, 2001.

\bibitem{zoulibiricci}
Yao~Wen Li and Xiao~Rong Zou.
\newblock {On the bi-Ricci curvature and some applications}.
\newblock {\em Houston journal of mathematics}, 34:467--481, 2008.

\bibitem{zbMATH06505276}
Hezi Lin.
\newblock On the structure of submanifolds in {Euclidean} space with flat
  normal bundle.
\newblock {\em Result. Math.}, 68(3-4):313--329, 2015.

\bibitem{lindqvist2017notes}
Peter Lindqvist.
\newblock {\em Notes on the p-Laplace equation}.
\newblock Number 161. University of Jyv{\"a}skyl{\"a}, 2017.

\bibitem{4}
Yozo Matsushima.
\newblock Vector bundle valued harmonic forms and immersions of {R}iemannian
  manifolds.
\newblock {\em Osaka Journal of Mathematics}, 8(1):1--13, 1971.

\bibitem{zbMATH01123720}
P.~Petersen and G.~Wei.
\newblock Relative volume comparison with integral curvature bounds.
\newblock {\em Geom. Funct. Anal.}, 7(6):1031--1045, 1997.

\bibitem{prs}
Stefano Pigola, Marco Rigoli, and Alberto~G Setti.
\newblock Constancy of $p$-harmonic maps of finite $q$-energy into
  non-positively curved manifolds.
\newblock {\em Mathematische Zeitschrift}, 258(2):347--362, 2008.

\bibitem{pigola2008constancy}
Stefano Pigola, Marco Rigoli, and Alberto~G Setti.
\newblock Constancy of p-harmonic maps of finite q-energy into non-positively
  curved manifolds.
\newblock {\em Mathematische Zeitschrift}, 258(2):347--362, 2008.

\bibitem{zbMATH07180899}
Keomkyo Seo and Gabjin Yun.
\newblock Liouville-type theorems for weighted {{\(p\)}}-harmonic 1-forms and
  weighted {{\(p\)}}-harmonic maps.
\newblock {\em Pac. J. Math.}, 305(1):291--310, 2020.

\bibitem{shen1996stable}
Ying Shen and Rugang Ye.
\newblock On stable minimal surfaces in manifolds of positive bi-ricci
  curvatures.
\newblock {\em Duke Mathematical Journal}, 85(1):109--116, 1996.

\bibitem{shen1997geometry}
Ying Shen and Rugang Ye.
\newblock On the geometry and topology of manifolds of positive bi-ricci
  curvature.
\newblock {\em arXiv preprint dg-ga/9708014}, 1997.

\bibitem{MR1145657}
Hiroshi Takeuchi.
\newblock Stability and {L}iouville theorems of {$p$}-harmonic maps.
\newblock {\em Japan. J. Math. (N.S.)}, 17(2):317--332, 1991.

\bibitem{tanno1996l2}
Shukichi Tanno.
\newblock L2 harmonic forms and stability of minimal hypersurfaces.
\newblock {\em Journal of the Mathematical Society of Japan}, 48(4):761--768,
  1996.

\bibitem{Tolksdorf1984}
Peter Tolksdorf.
\newblock Regularity for a more general class of quasilinear elliptic
  equations.
\newblock {\em Journal of Differential Equations}, 51(1):126--150, 1984.

\bibitem{zbMATH07178581}
Lili Wang and Guofang Wei.
\newblock Local {Sobolev} constant estimate for integral {Bakry}-{\'e}mery
  {Ricci} curvature.
\newblock {\em Pac. J. Math.}, 300(1):233--256, 2019.

\bibitem{wang2006harmonic}
Qiaoling Wang.
\newblock Harmonic maps and the topology of manifolds with positive spectrum
  and stable minimal hypersurfaces.
\newblock {\em Publicacions Matem{\`a}tiques}, pages 301--313, 2006.

\bibitem{wangql}
Qiaoling Wang.
\newblock Complete submanifolds in manifolds of partially non-negative
  curvature.
\newblock {\em Annals of Global Analysis and Geometry}, 37(2):113--124, 2010.

\bibitem{wang2011local}
Xiaodong Wang and Lei Zhang.
\newblock Local gradient estimate for p-harmonic functions on riemannian
  manifolds.
\newblock {\em Communications in analysis and geometry}, 19(4):759--772, 2011.

\bibitem{zbMATH07024088}
Jia-Yong Wu.
\newblock Comparison geometry for integral {Bakry}-{\'e}mery {Ricci} tensor
  bounds.
\newblock {\em J. Geom. Anal.}, 29(1):828--867, 2019.

\bibitem{yaust}
Shing~Tung Yau.
\newblock Some function-theoretic properties of complete {R}iemannian manifold
  and their applications to geometry.
\newblock {\em Indiana University Mathematics Journal}, 25(7):659--670, 1976.

\bibitem{Zhang2016A}
Jian~Feng Zhang and Yue Wang.
\newblock {A theorem of Liouville type for p-harmonic maps in weighted
  Riemannian manifolds}.
\newblock {\em Kodai Mathematical Journal}, 39(2):354--365, 2016.

\bibitem{9}
Xi~Zhang.
\newblock A note on $p$-harmonic 1-forms on complete manifolds.
\newblock {\em Canadian mathematical bulletin}, 44(3):376--384, 2001.

\end{thebibliography}
	\end{document}